\documentclass{amsart}

\usepackage{amsmath,amstext,amsthm,amsopn,amssymb,amscd}
\usepackage{tikz}
\usetikzlibrary{matrix}
\usetikzlibrary{arrows}
\usetikzlibrary{decorations.pathmorphing}

\usepackage{color}
\usepackage{pstricks, pstricks-add, pst-node}

 \newtheorem{thm}{Theorem}[section]
 \newtheorem{proposition}[thm]{Proposition}
 \newtheorem{corollary}[thm]{Corollary}
 \newtheorem{lemma}[thm]{Lemma}
 \newcommand{\dontprint}[1]{\relax}

\theoremstyle{definition}
\newtheorem{defi}[thm]{Definition}
 \newtheorem{cor}[thm]{Corollary}
\theoremstyle{remark}
 
 \newtheorem{rem}[thm]{Remark}

\newcommand{\alg}[1]{\mathfrak{{#1}}}
\newcommand{\tr}{\text{tr}}

\newcommand{\cn}[1]{\mathbb{C}^{{#1}}}

\newcommand{\co}[2]{\left[{#1},{#2}\right]} 

\newcommand{\morphU}{{\mathcal{U} }} 
\newcommand{\mU}{{\morphU }} 
\newcommand{\mV}{{\mathcal{V} }} 
\newcommand{\mW}{{\mathcal{W} }}

\newcommand{\Tpoly}{T_{\rm poly}}
\newcommand{\Dpoly}{D_{\rm poly}}
\newcommand{\formal}{\mathit{fml}}

\newcommand{\p}{\partial}

\newcommand{\C}{{\mathbb{C}}}
\newcommand{\R}{{\mathbb{R}}}

\DeclareMathOperator{\dv}{div}

\newcommand{\ndash}{\nobreakdash-\hspace{0pt}}

 \title{The character map in deformation quantization}
 \author{Alberto S. Cattaneo}
 \address{A. S. C.:
 Institut f\"ur Mathematik,
 Universit\"at Z\"urich-Irchel,
 Winterthurerstrasse 190,
 CH-8057 Z\"urich, Switzerland}
 \author{Giovanni Felder}
 \address{G. F.: Department of mathematics, ETH Zurich, CH-8092
 Zurich, Switzerland}
 \author{Thomas Willwacher}
 \address{T. W.: Department of Mathematics, Harvard University, Cambridge, MA}
 \thanks{This work been partially supported by SNF Grants 200020-121640/1 and 200020-122126/1, by the European
Union through the FP6 Marie Curie RTN ENIGMA (contract number
MRTN-CT-2004-5652), and by the European Science Foundation through
the MISGAM program.}

\begin{document}

\begin{abstract}
The third author recently proved that the Shoikhet--Dolgushev
$L_\infty$-morphism from Hochschild chains of the algebra of smooth
functions on a manifold to differential forms extends to cyclic
chains. Localization at a solution of the Maurer--Cartan equation
gives an isomorphism, which we call character map, from the periodic
cyclic homology of a formal associative deformation of the algebra
of functions to de Rham cohomology. We prove that the character map
is compatible with the Gauss--Manin connection, extending a result
of Calaque and Rossi on the compatibility with the cap product. As a
consequence, the image of the periodic cyclic cycle 1 is independent
of the deformation parameter and we compute it to be the A-roof
genus of the manifold.
Our results also imply the Tamarkin--Tsygan index Theorem.

\end{abstract}

\maketitle

\section{Introduction}
One of the consequences of the formality theorem for cyclic chains
\cite{Willwacher, Dolgushevetal} is the existence of a {\em
character map} from the periodic cyclic homology of any formal
associative deformation $A_\hbar$ of algebras of functions on a
manifold $M$ to the de Rham cohomology $H^\bullet(M,\mathbb R[[\hbar]])$.
By Kontsevich's formality theorem for Hochschild cochains such
deformations are classified by Poisson bivector fields in
$\hbar\Gamma(M,\wedge^2TM)[[\hbar]]$. A priori the character map
depends on the choice of formality map for Hochschild cochains and
cyclic chains. We consider here the Kontsevich $L_\infty$-morphism
of cochains and the Shoikhet--Dolgushev $L_\infty$-morphism of
modules from Hochschild chains to differential forms. The latter was
shown by the third author \cite{Willwacher} to be compatible with
the Rinehart--Connes differential $B$ and the de Rham differential
and therefore extends to a morphism of $L_\infty$-modules from
cyclic chains to the de Rham complex. In this paper we study the
dependence of the character map on the Poisson bivector field and
show that it is compatible with the Gauss--Manin connection. In
particular, we show that the image of the class of the cyclic cycle $1$ is
independent of the Poisson bivector field and can thus be computed
at the zero Poisson bracket, where it is given by the $A$-roof genus
of $M$. This proves one part of Tsygan's conjecture \cite{Tsygan} on
the character map (the independence on parameters) and disproves the
other (about characteristic classes of foliations) for the
Shoikhet--Dolgushev choice of formality map. Our main technical
result is that the formality map on periodic chains is compatible
with the Gauss--Manin connection. This result is an extension of the
Calaque--Rossi theorem \cite{CalaqueRossi} on compatibility of the
formality map on Hochschild chains with the cap product.

In the rest of the Introduction we recall the main notions (in
\ref{ss-1}) and describe the character map (in \ref{ss-2}). We then
state our main result on the compatibility of the character map with
the Gauss--Manin connection (Theorem \ref{t-1}) and give the formula
for the image of the cycle $1$ (Corollary \ref{cor:imof1}). The
proof of the independence of the image of $1$ on parameters is
easier in the case of regular Poisson structures; we give it in
Section \ref{s-2}. The evaluation of the image of $1$ at the zero
Poisson structure is discussed in Section \ref{s-3}. The case of
general Poisson structures is discussed in Section \ref{s-4}, 
where we also apply our results to obtain a proof of the Tamarkin--Tsygan
index Theorem.


\subsection*{Acknowledgements}
We are very grateful to Vasiliy Dolgushev and Boris Tsygan for helpful comments on the manuscript. They also independently obtained results similar to ours. 

\subsection{Cyclic chains as a module over Hochschild cochains and the Gauss--Manin
connection}\label{ss-1}
 Here we recall some basic notions and fix
sign conventions.

Let $C^\bullet(A)=\mathrm{Hom}_k(\bar A^{\otimes\bullet},A)$ be the
normalized Hochschild cochain complex of a unital associative
algebra $A$ over $k=\mathbb R$ or $\mathbb R[[\hbar]]$. Here $\bar
A=A/k1$. The degree-shifted complex $C^{\bullet+1}(A)$ with the
Gerstenhaber bracket
 \[
 [\phi,\psi]=\phi\bullet\psi-(-1)^{pq}\psi\bullet\phi,\qquad \phi\in
 C^{p+1}(A),\qquad\psi\in C^{q+1}(A),
 \]
 \[
 \phi\bullet\psi (a_1,\dots,a_{p+q+1})
 =\sum_{i=1}^{p+1}
 (-1)^{(i-1)q}\phi(a_1,\dots,a_{i-1},\psi(a_{i},\dots,a_{i+q}),\dots,a_{p+q+1}),
\]
is a differential graded Lie algebra; the differential is the
bracket $b=[\mu,\bullet]$ with the product $\mu\in C^2(A)$. The Lie
bracket induces a Lie bracket on the cohomology
$\mathit{HH}^{\bullet+1}(A)$ of $C^{\bullet+1}(A)$.

The complex of normalized Hochschild chains $C_\bullet(A)=A\otimes
\bar A^{\otimes(-\bullet)}$ is a graded module over the dgla
$C^\bullet(A)$. It is concentrated in non-positive degrees. The
action of a cochain $\phi\in C^{p+1}(A)$ on a chain
$a=(a_0,\dots,a_n)\in C_{-n}(A)$ is
\begin{eqnarray*}
L_\phi a&=&
(-1)^p
 \sum_{i=0}^{n-p}
(-1)^{ip}(a_0,\dots,a_{i-1},\phi(a_i,\dots,a_{i+p}),\dots,a_n)\\
&&+
(-1)^p
 \sum_{i=n-p+1}^n(-1)^{i n}
(\phi(a_i,\dots,a_n,a_0,\dots,a_{p-n+i-1}),\dots,a_{i-1}).
\end{eqnarray*}
The Hochschild differential on $C_\bullet(A)$ is the action of the
product $b=L_\mu$. The module property implies that $b\circ
L_\phi-(-1)^p L_\phi\circ b=L_{b\phi}$ for all $\phi\in C^{p+1}(A)$,
so that the action induces an action of the graded Lie algebra
$\mathit{HH}^{\bullet+1}(A)$ on the homology
$\mathit{HH}_{\bullet}(A)$ of $(C_\bullet,b)$. There is a second
differential $B$ of degree $-1$ on this complex (anti-)commuting
with $b$: the Rinehart--Connes differential
\[
B(a_0,\dots,a_n)=\sum_{i=0}^n(-1)^{i n}
(1,a_i,\dots,a_n,a_0,\dots,a_{i-1}).
\]
The induced action of $\phi\in \mathit{HH}^{p}(A)$ on
$\mathit{HH}_\bullet(A)$ obeys the homotopy formula $L_\phi=
B\circ I_\phi-(-1)^pI_\phi\circ B$, where $I_\phi\colon C_{-n}(A)\to
C_{-n+p}(A)$ is the internal multiplication, defined on chains by
\begin{equation}
\label{equ:Idefhoch}
I_\phi(a_0,\dots,a_{n})=(a_0\phi(a_1,\dots,a_p),a_{p+1},\dots,a_n).
\end{equation}
Here is a version of the homotopy formula on chains due to Getzler
\cite{Getzler}, who extended Rinehart's formula \cite{Rinehart} for
$p=1$.
\begin{lemma}\label{lemma1}
Let $\phi\in C^p(A)$ and let $H_\phi\colon C_{-n}(A)\to
C_{-n+p-2}(A)$ be the map
\begin{eqnarray*}
H_\phi(a_0,\dots,a_n)&=&\sum_{i=0}^{n-p}\sum_{j=0}^{n-p-i}
(-1)^{i(n-p+1)+j(p+1)}\\ &&
(1,a_{n+1-i},\dots,a_n,a_0,\dots,a_j,\phi(a_{j+1},\dots,a_{j+p}),\dots,a_{n-i})
\end{eqnarray*}
Then $H_\phi$ commutes with $B$ and
\[
 L_\phi=[B,I_\phi]-H_{b\phi}+[b,H_\phi].
\]
Here $[\ ,\ ]$ denotes the graded commutator in
$\mathrm{End}_k(C_\bullet(A))$.
\end{lemma}

\begin{proof}
Both $H_\phi\circ B$ and $B\circ H_\phi$ vanish since
$(1,1,\dots)=0$ in the normalized complex. The formula for $L_\phi$
can be proved by explicit calculations, see \cite{Rinehart,Getzler}.
\end{proof}

 The periodic cyclic chain complex is the complex of Laurent
series $PC_\bullet(A)=C_\bullet(A)((u))$ in an
indeterminate $u$ of degree 2 and $u$-linear differential $D=b+uB$.
Lemma \ref{lemma1} and the identities $[b,I_\phi]=I_{b\phi}$,
$[B,H_\phi]=0$ imply that the $u$-linear extension of $L_\phi$ obeys
\begin{equation}\label{e-G}
L_\phi=\frac 1 u ([D,\hat I_\phi]-\hat I_{b\phi}),\qquad
\hat I_\phi=I_\phi+u H_\phi,
\end{equation}
on the periodic cyclic complex.

Let $(A_t)$ be a one-parameter family of algebras given by an
algebra over $k[t]$ equal to $A[t]$ as a $k[t]$-module. Then the
product has the form $\mu+\gamma(t)$, where $\mu$ is the product in
$A=A_0$ and $\gamma(t)$ is a solution of the Maurer--Cartan equation
$b\gamma(t)+\frac12[\gamma(t),\gamma(t)]=0$. As noticed by Getzler
\cite{Getzler}, eq.~\eqref{e-G} implies that there exists a
connection, called Gauss--Manin connection, on the periodic cyclic
cohomology viewed as a vector bundle over the parameter space. The
parallel translation is given at the level of chains by the
following formula.
\begin{lemma}\label{l-1}{\rm(Getzler)}
Suppose that $c(t)\in PC_\bullet(A_t)$ is a solution of the
differential equation
\[
\frac d{dt}c(t){}+\frac 1u\hat I_{\dot\gamma(t)}c(t)=0,
\]
and that $c(0)$ is a cyclic cycle. Then $c(t)$ is a cyclic cycle for
all $t$ and its class in $PH_\bullet(A_t)$ depends only on  the
class of $c(0)$.
\end{lemma}

The statement follows from \eqref{e-G} and the fact that
$\dot\gamma(t)=\frac d{dt}\gamma(t)$ is a 2-cocycle in $C^2(A_t)$.

Getzler also showed that this formula defines a flat connection on
the periodic cyclic cohomology of any family of algebras. Tsygan
\cite{TsyganGM} extended the Gauss--Manin connection to a flat
superconnection at the level of chains, see also
\cite{Dolgushevetal2009}.

\subsection{The character map}\label{ss-2}
Let now $k=\mathbb R$ and $A=C^\infty(M)$ be the algebra of smooth
functions on a manifold $M$ and $C^\bullet(A)$ be the complex of
multidifferential Hochschild cochains. All operations described
above are defined on this subcomplex of the Hochschild complex.

Kontsevich \cite{Kontsevich}, in his proof of the formality theorem,
constructed an $L_\infty$-quasi\-isomorphism $\mathcal U$ from the
dgla $T^{\bullet+1}(M)=\Gamma(M,\wedge^{\bullet+1} TM)\simeq
\mathit{HH}^{\bullet+1}(A)$, with Schouten--Nijenhuis
bracket\footnote{The Schouten--Nijenhuis bracket is the Lie bracket
on vector fields and is extended to general multivector fields by
the rule
$[\alpha\wedge\beta,\gamma]=\alpha\wedge[\beta,\gamma]+(-1)^{pq}\beta\wedge[\alpha,\gamma]$,
$\alpha\in T^p(M),\beta\in T^q(M)$. With this sign convention the
HKR homomorphism sends $\xi_1\wedge\cdots\wedge\xi_p$ to the
Hochschild cocycle
$a_1\!\otimes\!\cdots\!\otimes\!a_p\mapsto(1/p!)\sum_{\sigma\in
S_p}\mathrm{sgn}(\sigma)\xi_{\sigma(p)}(a_1)\cdots\xi_{\sigma(1)}(a_p)$}
and zero differential, to $C^{\bullet+1}(A)$. The first component
$\mathcal U_1$ is the Hochschild--Kostant--Rosenberg (HKR)
quasi-isomorphism. Shoikhet \cite{Shoikhet} proved Tsygan's
formality conjecture for chains on $M=\mathbb R^n$ by giving a
morphism $\mathcal V$ of $L_\infty$-modules over $T^{\bullet+1}(M)$
from $C_\bullet(A)$ to the dg module of differential forms
$\Omega^{-\bullet}(M,\mathbb R)$ with zero differential. The action
of $\gamma\in T^{p}(M)$ on $\Omega^{-\bullet}$ is given by the Lie
derivative
$L_\gamma=d\circ\iota_\gamma-(-1)^p\iota_\gamma\circ d$, where $d$
is the de Rham differential and the internal multiplication of
vector fields is extended to multivector fields by the rule
$\iota_\gamma\iota_\eta=\iota_{\gamma\wedge\eta}$. The
$L_\infty$-action of $T^{\bullet+1}(M)$ on Hochschild chains is the
pull-back of the action of cochains on chains by Kontsevich's
morphism $\mathcal U$. Shoikhet's result was extended to general
manifolds by Dolgushev \cite{Dolgushev}. The first component
$\mathcal V_0$ of the $L_\infty$-morphism $\mathcal V$ is the HKR
map $(a_0,\dots,a_n)\mapsto (1/n!)a_0da_1\cdots da_n$ and induces an
isomorphism in homology if we interpret $C_{-n}(A)$ as the complex
of $\infty$-jets around the small diagonal of $M^{n+1}$. Finally, it
was shown by the third author \cite{Willwacher} that
Shoikhet's (local) morphism commutes with the Connes differential
$B$ on chains and the de Rham differential on differential forms.
Globalizing, one obtains a morphism of $L_\infty$-modules over
$T^{\bullet+1}(M)$ from the periodic cyclic complex
$PC_{\bullet}(A)$ to the de Rham complex
$(\Omega^{-\bullet}(M,\mathbb R[[\hbar]])((u)),ud)$. We will henceforth 
denote this morphism by $\mV$, and refer to it as the cyclic Shoikhet-Dolgushev morphism.
We advise the reader that this is slightly misleading. While in the local 
case our $\mV$ is indeed identical to Shoikhet's morphism, the globalized version 
differs since the globalization procedure is sensitive to the change in the
differential from $C_{-\bullet}(A)$ to $PC_{\bullet}(A)$.

For any solution of the Maurer--Cartan equation $[\pi,\pi]=0$ in
$\hbar T^2(M)[[\hbar]]$, i.e., a formal Poisson bivector field, we
get a solution $\mathcal U_0^\pi=\sum_{n=1}^\infty \frac1{n!}\mathcal
U_n(\pi,\dots,\pi)$ of the Maurer--Cartan equation
$b\,\mathcal{U}_0^{\pi}+\frac12[\mathcal{U}_0^{\pi},\mathcal{U}_0^{\pi}]=0$
in $C^2(A)[[\hbar]]$. The Maurer--Cartan equation is the
associativity of the star product $f\star
g=fg+\mathcal{U}_0^\pi(f\otimes g)$ on $A[[\hbar]]$. Let $A_\hbar$
be the algebra $(A[[\hbar]],\star)$. The twists at $\pi$ of
$\mathcal U,\mathcal V$ are $L_\infty$-morphisms $\mathcal
U^\pi,\mathcal V^\pi$ defined by
\[
 \mathcal U_n^\pi(\gamma_1,\dots,\gamma_n)= \sum_{m=0}^\infty
 \frac{1}{m!}
 \mathcal U_{n+m}(\gamma_1,\dots,\gamma_n,\pi,\dots,\pi),\qquad n\geq 1,
\]
\[
 \mathcal V_n^\pi(\gamma_1,\dots,\gamma_n)= \sum_{m=0}^\infty
 \frac{1}{m!}
 \mathcal V_{n+m}(\gamma_1,\dots,\gamma_n,\pi,\dots,\pi),\qquad n\geq0.
\]
The twist $\mathcal U^\pi$ is an $\hbar$-linear
$L_\infty$-quasiisomorphism from $T^{\bullet+1}(M)[[\hbar]]$ with
Poisson differential $[\pi,\ ]$ to $C^{\bullet+1}(M)[[\hbar]]$ with
Hochschild differential $b_\star=b+L_{\mathcal{U}_0^\pi}$ calculated
with the star product. Similarly $\mathcal V^\pi$ is an
$L_\infty$-morphism of modules over the dgla
$(T^{\bullet+1}(M),[\pi,\ ])$ from $(PC_{\bullet}(A)[[\hbar]],
b_\star+uB)$ to $(\Omega^{-\bullet}(M,\mathbb
R[[\hbar]])((u)),L_\pi+u d)$.
 We refer to \cite{Yekutieli, Kontsevich, Tsygan, Shoikhet, Dolgushev,
Willwacher} for more details.

Since
$L_\pi=d\circ\iota_\pi-\iota_\pi\circ d$, the map $\alpha \mapsto
e^{\iota_\pi/u}\alpha$ %
 is an isomorphism of complexes
\[
 (\Omega^{-\bullet}(M,\mathbb R[[\hbar]])((u)),L_\pi+ud)
 \to
 (\Omega^{-\bullet}(M,\mathbb R[[\hbar]])((u)),u d)
\]
Passing to cohomology, we obtain the {\em character map}
\[
\mathit{PH}_{\bullet}(A_\hbar)\to H^{-\bullet}(M,\mathbb
R[[\hbar]])((u)),
\]
of modules over the Poisson cohomology, induced by
\[
\tilde{\mathcal V}^\pi_0=e^{\iota_\pi/u}\circ \mathcal V^\pi_0, 
\]
from the periodic cyclic cohomology of the deformed algebra
$A_\hbar$ to de Rham cohomology.

\subsection{Dependence on parameters}\label{ss-3}
Suppose we have a one-parameter family of solutions $\pi(t)$ of the
Maurer--Cartan equation, with, say, polynomial dependence on $t$,
for example $\pi(t)=t\pi$ for some Poisson structure $\pi$. We then
obtain a one-parameter family of star-products whose reduction
modulo $\hbar^N$ depends polynomially on $t$ for any $N$. The
Gauss--Manin connection allows us to identify the periodic cyclic
cohomology for different values of the parameter, so we can study
the dependence on $t$ of the character map.

\begin{thm}\label{t-1}
 Let $A_t$ be the quantum algebra of functions on $M$
associated with the Poisson bracket $\pi(t)$. Suppose $c(t)$ is a
cycle in $PC_{-n}(A_t)$, horizontal with respect to the Gauss--Manin
connection, i.e., such that
\[
\frac d{dt}c(t){}+\frac1u\hat
I_{\mathcal{U}_1^{\pi(t)}(\dot\pi(t))}c(t)=0,\qquad \hat
I_\gamma=I_\gamma+u H_\gamma.
\]
Then the class of $\tilde{\mathcal V}^{\pi(t)}_0\left(c(t)\right)$
in $\oplus_mH^{n+2m}(M,\mathbb R[[\hbar]])u^m$ is independent of
$t$.
\end{thm}

 The proof of this theorem is based on an identity extending the
``compatibility with cap product'' result of Calaque and Rossi
\cite{CalaqueRossi}. To formulate the result we introduce the
modified internal multiplication of a multivector field $\gamma$ on
differential forms:
\[
\hat\iota_\gamma=\iota_\gamma+u\frac12dL_\gamma.
\]
It is clear that the Cartan formula $L_\gamma=[d, \hat\iota_\gamma]$ still holds. %

\begin{proposition}\label{prop2}
Let $\pi$ be a Poisson bivector field on $M$, $A=C^\infty(M)$ and
$\gamma\in \Gamma(M,\wedge^pTM)$ such that $[\pi,\gamma]=0$.
\begin{enumerate}
 \item The map 
\[
 \hat\iota_\gamma\circ\mathcal V_0^\pi
 -\mathcal V_0^\pi\circ\hat I_{\mathcal U_1^{\pi}(\gamma)}
 +u\mathcal V_1^\pi(\gamma)
 \]
is a map of complexes from the periodic cyclic complex $PC_\bullet(A_\hbar)$ (with differential $b_\star+uB$) to the complex $\Omega^{-\bullet}((u))[[\hbar]]$ (with differential $L_\pi + \rm u \, d$). Furthermore this map induces the zero map on cohomology.
 \item In the special case $M=\R^n$ and $\pi, \gamma$ as before, 
 there exist maps
\[
 \mathcal V_I^\pi(\gamma)\colon C_\bullet(A)\to \Omega^{-\bullet+p-1}
 (M,\mathbb R[[\hbar]]),\quad
 \mathcal V_H^\pi(\gamma)\colon C_\bullet(A)\to \Omega^{-\bullet+p-3}
 (M,\mathbb R[[\hbar]]),
 \]
  linear in $\gamma$, such that
\[
 \hat\iota_\gamma\circ\mathcal V_0^\pi
 -\mathcal V_0^\pi\circ\hat I_{\mathcal U_1^{\pi}(\gamma)}
 +u\mathcal V_1^\pi(\gamma)
 =(L_\pi+ud)\circ
 X^\pi(\gamma)+X^\pi(\gamma)\circ(b_\star+uB),
 \]
 with homotopy
$X^\pi(\gamma)=\mathcal V_I^\pi(\gamma)+u\mathcal V_H^\pi(\gamma)$.
\end{enumerate}
\end{proposition}
\begin{proof}[Simple part of the proof]
 Let us compute the commutator with the differential
\begin{multline*}
 (b_\star+uB)\circ (\hat\iota_\gamma\circ\mathcal V_0^\pi
 -\mathcal V_0^\pi\circ\hat I_{\mathcal U_1^{\pi}(\gamma)}
 +u\mathcal V_1^\pi(\gamma) ) 
-(-1)^p
(\hat\iota_\gamma\circ\mathcal V_0^\pi
 -\mathcal V_0^\pi\circ\hat I_{\mathcal U_1^{\pi}(\gamma)}
 +u\mathcal V_1^\pi(\gamma) ) \circ (u\rm d)
\\=u L_\gamma \circ\mathcal V_0^\pi
- u \mathcal V_0^\pi \circ L_\gamma
 + u (V_0^\pi\circ L_\gamma - V_0^\pi\circ L_\gamma) = 0.
\end{multline*}
\end{proof}

%

Proposition \ref{prop2} is our main technical result and is proven
in Section \ref{s-4}. The operators $\mathcal V_I$ and $\mathcal
V_H$ are given explicitly in terms of integrals over configuration
spaces. Given this proposition, the proof of the theorem is a simple
calculation. With the notation  $\pi=\pi(t), \dot\pi=\frac
d{dt}\pi(t), b_\star(t)=b+L_{\pi(t)}$ we have:
\begin{eqnarray*}
 \frac d{dt}{\mathcal V}_0^{\pi}
 +\frac1u\hat\iota_{\dot\pi}\circ{\mathcal V}_0^{\pi}
 &=&
 \mathcal V_1^{\pi}(\dot\pi){}+\frac1u
\hat\iota_{\dot\pi}\circ{\mathcal V}_0^{\pi}
 \\
 &=&
 \frac1u
 \bigl( \mathcal V_0^{\pi}\circ\hat I_{\mathcal U_1^{\pi}(\dot
 \pi)}
  +
 (L_{\pi}+u d)\circ X^{\pi}(\dot\pi)
 \\
 &&+
 X^\pi(\dot\pi)\circ (b_\star(t)+u B)
 \bigr)
\end{eqnarray*}
If $c(t)$ is a cycle obeying the differential equation we get
\[
 \left(\frac d{dt}{}+\frac1u\hat\iota_{\dot\pi}\right)\left({\mathcal V}_0^{\pi}c(t)\right)
 =
 \left(d+\frac1uL_\pi\right)\left(X^{\pi}(\dot\pi) c(t)\right)
\]
This formula can be rewritten as
\[
 \frac d{dt}\left(e^{\hat\iota_{\pi}/u}{\mathcal V}_0^{\pi}c(t)\right)
 =
 d\left(e^{\hat\iota_{\pi}/u}X^{\pi}(\dot\pi) c(t)\right).
\]
Thus the de Rham class of $e^{\hat\iota_{\pi}/u}{\mathcal
V}_0^{\pi}c(t)$ is independent of $t$. Moreover, since
$e^{\hat\iota_{\pi}/u}=(1-\frac1{2u} d\iota_{\pi}
d)e^{\iota_{\pi}/u}$, and $e^{\iota_{\pi}/u}{\mathcal
V}_0^{\pi}c(t)$ is a cocycle, we may replace $\hat\iota_{\pi}$ by
$\iota_\pi$ in the exponential.

\noindent{\bf Remark.} The identity of Proposition \ref{prop2}
reduces at $u=0$ to the identity
\[
 \iota_\gamma\circ\mathcal V_0^\pi
 -\mathcal V_0^\pi\circ I_{\mathcal U_1^{\pi}(\gamma)}
 =
 L_\pi\circ\mathcal V_I^\pi(\gamma)
 +\mathcal V_I^\pi(\gamma)\circ b_\star.
\]
This is the statement of compatibility with the cap product
$(I_\phi=\phi\,\cap)$ proved in \cite{CalaqueRossi}.

\medskip

\noindent{\bf Remark.} Geometrically one should think of the
periodic cyclic complex $PC_\bullet(A_\hbar)$ as the fibre of a
trivial vector bundle on the Maurer--Cartan variety of formal
deformations of the product. This vector bundle carries the
Gauss--Manin connection $d{}+\theta$, where $\theta$ is a one-form
with values in the endomorphisms of $C_\bullet$. The value of this
one-form on a tangent vector $\gamma\in \mathrm{Ker}(b_\star)$ is
\[
\theta_\gamma=\frac1u\hat I_\gamma=\frac1u I_\gamma+H_\gamma.
\]
The statement of Proposition \ref{prop2} is that the
Shoikhet--Dolgushev quasi-isomorphism intertwines, up to homotopy,
this connection and the connection $d{}+\frac1u\hat\iota$ on the
trivial vector bundle on the Maurer--Cartan variety of Poisson
structure with fibre $\Omega^{-\bullet}(M,\mathbb
R[[\hbar]])((u))$.

\subsection{The image of the cycle 1}\label{ss-4}
\label{sec:imof1}
For any unital algebra $A$, $1$ is a cyclic cycle and thus also a
periodic cyclic cycle.
\begin{corollary}
\label{cor:imof1} Let $\pi\in\hbar \Gamma(M,\wedge^2TM)[[\hbar]]$ be
a formal Poisson bivector field. The image of the cyclic cycle $1\in
(C^\infty(M)[[\hbar]],\star)$ by the character map is $\hat
A_0(M)+u\hat A_2(M)+u^2\hat A_4(M)+\cdots$, where $\hat A_i(M)\in
H^i(M)$ are the components of the $A$-roof genus of $M$.
\end{corollary}

Indeed, $\pi(t)=t\pi$ is a Poisson bivector field for any $t$. By
Theorem \ref{t-1} we can evaluate the class of $\tilde {\mathcal
V}_0^{t\pi}(1)$ at $t=0$. We do this in Section \ref{s-3}.

%
%

\section{Special case: Regular Poisson structures}
\label{s-2}
In this section, we will illustrate the above results for regular Poisson structures. In this special case, there is a simple proof of $\hbar$-independence by elementary means that does not require the more intricate Theorem \ref{thm:compat}. The results of this section are not needed anywhere else in the paper.

\begin{proposition}
\label{prop:regularcase}
For $\mV$ the cyclic Shoikhet--Dolgushev morphism and $\pi$ a regular Poisson structure, the de Rham cohomology class of
\[
 \tilde{\mV}_0(1) = e^{\hbar \iota_\pi / u} \mV_0^\pi(1)
\]
is independent of $\hbar$.
\end{proposition}

In general, let us introduce the notation $\tilde{\mV}_j = e^{\hbar\iota_\pi/u}\mV_j^\pi$. It is then easy to see that the first $L_\infty$ relations reduce to the following 
\begin{align*}
u{\rm d} \tilde{\mV}_0(\alpha) &= \tilde{\mV}_0((b_\star+uB)\alpha)\\
 u{\rm d} \tilde{\mV}_1(\gamma;\alpha)+ (-1)^{|\gamma|}(L_\gamma+\hbar \iota_{\co{\pi}{\gamma}};\alpha) \tilde{\mV}_0(\alpha)
 &=
 (-1)^{|\gamma|} \tilde{\mV}_1(\gamma;(b_\star+uB)\alpha) 
\\
&\quad +\hbar \tilde{\mV}_1(\co{\pi}{\gamma};\alpha) + (-1)^{|\gamma|}\tilde{\mV}_0( L_{\mU(\gamma)}\alpha)
\end{align*}

For the proof of the proposition, first note that the subcomplex $(T^\bullet_\parallel,\co{\pi}{\cdot}) \subset (T^\bullet,\co{\pi}{\cdot})$ of multivector fields tangential to the symplectic foliation is locally acyclic. Concretely, in a small enough neighbourhood around any point, the manifold looks like a product of a symplectic space and one with trivial Poisson structure.

It follows that there is a good covering $\{U_i\}_{i\in I}$ of $M$, a collection of (locally defined) tangential vector fields $\xi=\{\xi_i\}_{i\in I}$, and a collection of (locally defined) functions $f=\{f_{ij}\}_{i\neq j\in I}$, $f_{ij}=-f_{ji}$, such that 
\begin{align*}
\pi &= \co{\pi}{\xi_i} &\text{on $U_i$} \\
\xi_i-\xi_j &= \co{\pi}{f_{ij}} &\text{on $U_i\cap U_j$} \, .
\end{align*} 

In other words, the class $[\pi]$ in Poisson cohomology can be represented by the cocycle $\check{\delta} f$ in the Poisson-\v Cech double complex, where $\check{\delta}$ is the \v Cech-differential.

Next compute on $U_i$:\footnote{Implicitly we use here that $\tilde \mV$ is local, i.e., its value at a point depends only on the jets of its arguments at that point. Hence it makes sense to restrict $\tilde \mV$ to $U_i$.}
\begin{align*}
\frac{d}{d\hbar} \tilde{\mV}_0(1)
&=
\frac{d}{d\hbar}
e^{\hbar \iota_\pi/u}
\sum_{j\geq 0}
\frac{\hbar^j}{j!}
\mV_j(\pi, \dots, \pi; 1)
\\
&=
\iota_\pi \tilde{\mV}_0(1)/u
+ \tilde{\mV}_1(\pi; 1)
\\
&=
\iota_\pi \tilde{\mV}_0(1)/u
+ \tilde{\mV}_1(\co{\pi}{\xi_i}; 1)
\\
&=
\iota_\pi \tilde{\mV}_0(1)/u
+
(-\frac{1}{\hbar} L_{\xi_i} -\iota_\pi/u) \tilde{\mV}_0(1)
+\frac{u}{\hbar} {\rm d} \tilde{\mV}_1(\xi_i;1)
\\
&=
\frac{1}{\hbar}
{\rm d} \left(-\iota_{\xi_i} \tilde{\mV}_0(1)
+u\tilde{\mV}_1(\xi_i;1) \right)
\end{align*}
The class of this cocycle in \v Cech-de Rham cohomology is the same as that of the \v Cech-degree 1 chain (denoting $\xi_{ij}=\xi_i-\xi_j$) 
\begin{align*}
\frac{1}{\hbar}
\left(
-\iota_{\xi_{ij}} \tilde{\mV}_0(1)
+ u \tilde{\mV}_1(\xi_{ij};1)
\right)
&=
\frac{1}{\hbar}
\left(
-\iota_{\xi_{ij}} \tilde{\mV}_0(1)
+u \tilde{\mV}_1(\co{\pi}{f_{ij}};1)
\right)
\\
&=
\frac{u}{\hbar^2}
\left(
{\rm d} \iota_{f_{ij}}\tilde{\mV}_0(1)
+
u {\rm d} \tilde{\mV}_1(f_{ij};1)
-
\tilde{\mV}_0(L_{\mU_1(f_{ij})} 1)
\right)
\\
&=
\frac{u}{\hbar^2}
\left(
{\rm d} \iota_{f_{ij}}\tilde{\mV}_0(1)
+
u {\rm d} \tilde{\mV}_1(f_{ij};1)
\phantom{\frac{1}{u}} \right.
\\ &\quad\quad\quad\quad\quad\quad\quad \left.
- \frac{1}{u}\tilde{\mV}_0((uB+b_\star) \mU_1(f_{ij}))
\right)
\\
&=
\frac{u}{\hbar^2}
{\rm d}
\left(
\iota_{f_{ij}}\tilde{\mV}_0(1)
+
u\tilde{\mV}_1(f_{ij};1)
-
\tilde{\mV}_0(\mU_1(f_{ij}))
\right)
\end{align*}
This cochain in turn describes the same \v Cech-de Rham cohomology class as the \v Cech-degree 2 cochain (denoting $f_{ijk}=f_{ij}+f_{jk}+f_{ki}$)
\begin{align*}
\frac{u}{\hbar^2}
\left( 
f_{ijk} \tilde{\mV}_0(1) 
+
u\tilde{\mV}_1(f_{ijk};1)
-
\tilde{\mV}_0( \mU_1(f_{ijk}))
\right)
\end{align*}
We claim that this expression is zero, at least for suitable choices made in the globalization step of the construction of $\mU$, $\mV$.\footnote{Note that, due to non-canonical choices to be made in the globalization step (see section \ref{sec:globreminder}), $\mU$ and $\mV$ are defined uniquely only up to homotopy.} Concretely, it follows from the following lemma.

\begin{lemma}
Let $\pi$ be a regular Poisson structure, $f$ be a Poisson central function, and $a_0,\dots,a_n$ be arbitrary functions. Let $\mU$ be the Kontsevich formality morphism and $\mV$ the cyclic Shoikhet-Dolgushev morphism, and $\mU^\pi$ and $\mV^\pi$ their twisted versions as before. Then, for suitable choices taken in the globalization step of $\mU$, $\mV$ the following hold:
\begin{gather}
\mV^\pi_1(f;a_0,\dots, a_n)=0 \\
\mV^\pi_0(f)=f\mV^\pi_0(1, a_1,\dots, a_n) \\
\mU^\pi_1(f)=f\, .
\end{gather}
\end{lemma}
\begin{proof}[Sketch of proof]
For the proof, we assume familiarity with M. Kontsevich's construction of the formality morphism $\mU^\pi$ and the construction of $\mV^\pi$ by V. Dolgushev and B. Shoikhet, see \cite{Kontsevich, Dolgushev,Shoikhet}. Both constructions proceed by first defining the morphisms in the local cases (i.e., $M=\R^n$), and then globalizing these morphisms. The globalization step is briefly recalled in section \ref{sec:globreminder}. It depends on the choice of a particular Fedosov connection. We will choose this connection in such a manner that the flat lift of $\pi$, say $\hat \pi$, is fiberwise constant. This is possible only for regular Poisson structures. 

Let us start with the third assertion of the Lemma. Using the notation from section \ref{sec:globreminder} below, $\mU^\pi$ is defined as the twisted version of the composition of quasi-isomorphisms:
\[
 \begin{CD}
\Tpoly @>F>> \Tpoly^\formal @>>> \Dpoly^\formal @<G<< \Dpoly.
\end{CD}
\]
Here the right hand morphism must be inverted up to homotopy before composing. To see that $\mU^\pi_1(f)=f$ for a Poisson central function $f\in\Tpoly^0$ one needs to trace it through these three morphisms. The first morphism $F$ only has a single Taylor component, and maps $f$ to its flat lift $\hat f$. Furthermore $F(\pi)=\hat \pi$ is fiberwise constant by our choice of Fedosov connection. The second morphism in the composition is the fiberwise application of M. Kontsevich's (local) formality morphism. It is defined by a sum of terms labelled by graphs. Using that $\hat\pi$ is fiberwise constant, that $\hat \pi$ is a flat lift and that $\co{\hat\pi}{\hat f}=0$, it is an elementary exercise to show that all Kontsevich graphs except the trivial one yield zero contribution. Hence the image of $f$ under the twisted version of the first two maps is $\hat f$. Finally one has to invert $G$, up to homotopy. There are explicit formulas for this inverse. Twisting and inserting $\hat f$ one can check explicitly that $\hat f$ is mapped to $f$.

Next consider the statements involving $\mV$. It is defined through the zig-zag of quasi-isomorphisms of $L_\infty$ modules (see again section \ref{sec:globreminder} for the notation)
\[
\begin{CD}
\Omega((u)) @>>> \Omega^\formal((u)) @<<< PC^\formal @<<< PC.
\end{CD}
\]
The left hand morphism has to be inverted before composing. Let us start with the chain $f\otimes a_1\otimes\dots\otimes a_n$, as an element of the periodic cyclic chain complex $PC$ on the right. The twisted version of the right hand map sends it to the flat lift $\hat f \otimes \hat a_1\otimes\dots\otimes \hat a_n$, plus some correction terms. The correction terms arise because the $L_\infty$ morphism is not just the flat lift, but picks up higher components through the inversion of the map of $L_\infty$ algebras $\Dpoly\to \Dpoly^\formal$, see \ref{sec:globreminder}.
The correction terms however can be seen to all have the form $\hat f \otimes (\cdots)$, i.e., the $\hat f$ is never hit by any derivatives.

The second morphism in the composition is defined using Shoikhet graphs. An example graph contributing to $\mV^\pi$ is shown in Figure \ref{fig:wheels}. However, note that all graphs in which an arrow hits the vertex $\bar 0$, corresponding to $\hat f$, contribute zero. This can be shown using that $f$ is locally constant along the symplectic leaves, and that $\hat \pi$ is fiber-wise constant. Hence the image of $f$ under the first two maps is the image of $1\otimes a_1\otimes\dots\otimes a_n$ times $\hat f$. Inverting the left hand map, one can check that it sends this product to $f\mV^\pi_0(1\otimes a_1\otimes\dots\otimes a_n)$ and the second statement of the lemma is shown. 

Finally let us turn to the first statement of the lemma. The component $\mV^\pi_1$ can receive contributions from each of the three morphisms in the composition. We claim that all these contributions are zero. To check this for the first and last morphism is technical, and will be left to the reader. Let us check that the second morphism contributes zero. Consider a (potentially) contributing Shoikhet graph, with, say $k\geq 1$ type I vertices, one of which is the one corresponding to $\hat f$. Vertices can have either one or two outgoing edges, furthermore there can be outgoing edges from the central vertex. However, using that $\hat \pi$ is fiberwise constant one sees that the maximal number of edges is one less than the dimension of the relevant configuration space, and hence the weight of the graph will be zero.
\end{proof}


\section{The cyclic Shoikhet--Dolgushev map for $\pi=0$}
\label{sec:hbar0}\label{s-3}
\begin{figure}
 \psset{unit=2cm}
\psset{arrowscale=1}

\begin{pspicture}(-2,-2)(2,2) 
  \cnode*(0,0){3pt}{G}
  \cnode*(2;  0){3pt}{a0}
  \cnode*(2;-30){3pt}{a1}
  \cnode*(2;-60){3pt}{a2}
  \cnode*(2;-90){3pt}{a3}
  \cnode*(2;-120){3pt}{a4}
  \cnode*[linecolor=gray](1;-75){3pt}{p1}
  \cnode*[linecolor=gray](1.5;-75){3pt}{p2}
  \cnode*[linecolor=gray](1.3;-10){3pt}{p3}
  \cnode*[linecolor=gray](1.4;-160){3pt}{p4}

  \pscircle(0,0){2}
  \psset{arrows=->}
  \ncline{G}{a1}
  \ncline{G}{a2}  
  \ncline{G}{a3}
  \ncline{G}{a4}
  \ncline{p1}{a3}
  \ncline{p1}{a2}
  \ncline{p2}{a3}
  \ncline{p2}{a2}
  \ncline{p3}{a0}
  \ncline{p3}{a1}
  \ncline{p4}{a4}  

  \uput{1ex}[0](a0){$\bar{0}$}
  \uput{1ex}[-90](a1){$\bar{4}$}
  \uput{1ex}[-90](a2){$\bar{3}$}
  \uput{1ex}[-90](a3){$\bar{2}$}
  \uput{1ex}[-90](a4){$\bar{1}$}

\rput(-.75,0.5){
  \cnode*(.5;30){3pt}{aw1}
  \cnode*(.5;-30){3pt}{aw2}
  \cnode*(.5;-90){3pt}{aw3}
  \cnode*(.5;-150){3pt}{aw4}
  \cnode*(.5;-210){3pt}{aw5}

  \rput[c](.5;90){\rnode{adots}{$\cdots$}}
}

  \psset{arrows=->}
  \ncline{aw1}{aw2}
  \ncline{aw2}{aw3}
  \ncline{aw3}{aw4}
  \ncline{aw4}{aw5}
  \ncline{aw5}{adots}
  \ncline{adots}{aw1}
  \ncline{G}{aw1}
  \ncline{G}{aw2}
  \ncline{G}{aw3}
  \ncline{G}{aw4}
  \ncline{G}{aw5}

\rput( .75,0.5){
  \cnode*(.5;30){3pt}{bw1}
  \cnode*(.5;-30){3pt}{bw2}
  \cnode*(.5;-90){3pt}{bw3}
  \cnode*(.5;-150){3pt}{bw4}
  \cnode*(.5;-210){3pt}{bw5}
  \rput[c](.5;90){\rnode{bdots}{$\cdots$}}
}

  \psset{arrows=->}
  \ncline{bw1}{bw2}
  \ncline{bw2}{bw3}
  \ncline{bw3}{bw4}
  \ncline{bw4}{bw5}
  \ncline{bw5}{bdots}
  \ncline{bdots}{bw1}
  \ncline{G}{bw1}
  \ncline{G}{bw2}
  \ncline{G}{bw3}
  \ncline{G}{bw4}
  \ncline{G}{bw5}
\ncline{p4}{aw4}

  \rput[c](0,1){\rnode{bigdots}{\huge$\,\, \cdots$}}


\end{pspicture}
\caption{\label{fig:wheels} The picture shown an example graph contributing to $\mV_0^\pi$ in the regular case. The grey vertices correspond to the (fiberwise constant) Poisson structure $\hbar \pi$. Thus putting $\hbar=0$ the only contributions come from graphs consisting entirely of the ``wheels''.}
\end{figure}
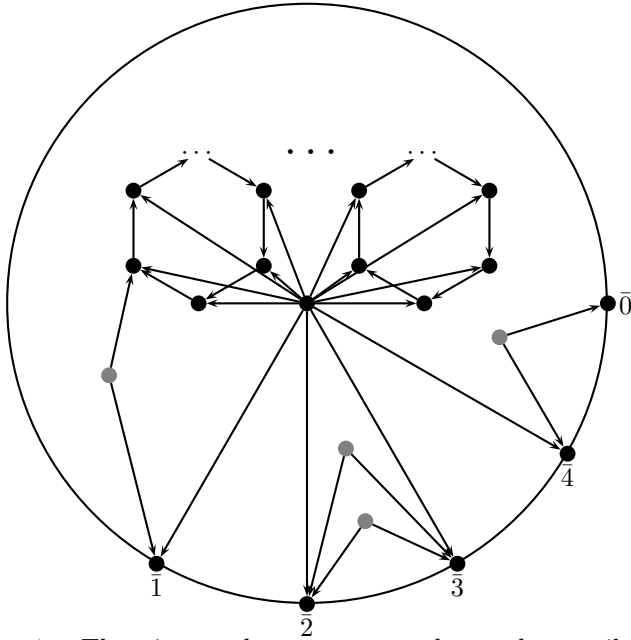

This section is devoted to the following result, needed in the proof of Corollary \ref{cor:imof1}.
\begin{proposition}
\label{prop:SDAhat}
For $\mV$ the cyclic Shoikhet--Dolgushev morphism, $\pi=0$, and $c$ a cyclic cycle in $PC_\bullet(A)$
\[
[\tilde{\mV}^{\pi=0}_0(c)] = \hat{A}_u(M)\, \mathit{Co}(c)\, .
\]
Here 
\[
\hat{A}_u(M):=\hat
A_0(M)+u\hat A_2(M)+u^2\hat A_4(M)+\cdots
\] with $\hat A_i(M)\in H^i(M)$ being the components of the $A$-roof genus of $M$. Furthermore $\mathit{Co}$ is the Connes map identifying periodic cyclic and de Rham cohomology, concretely on chains
\[
\mathit{Co}(a_0,\dots, a_n) = a_0da_1\wedge \dots \wedge da_n\, .
\]
\end{proposition}
\begin{proof}
The only contributing graphs are unions of wheels and edges directly connecting the center vertex with vertices on the boundary, as shown in Figure \ref{fig:wheels} (without the grey vertices). It is not hard to see that the latter (direct) edges produce the factor $\mathit{Co}(c)$. Let us turn to the wheels.  
Consider a single wheel with $j$ vertices. It has a number (weight) $w_j$ and a differential form $\alpha_j$ associated to it. The contribution of all combinations of wheels is hence 
\[
 [\tilde{\mV}_0(1)] = [e^{\sum_{j\geq 2}w_j \alpha_j }]\, .
\]
From the details of the globalization procedure (and tedious sign conventions) it can be seen that 
\[
 \alpha_j = (-1)^{j(j-1)/2}u^j\tr(R^j)  
\]
where $R$ is the curvature 2-form of some affine torsion free connection on $M$. Furthermore, the weights $w_j$ can be computed explicitly. A simple reflection argument shows that $w_j=0$ if $j$ is odd. For $j$ even it was shown by M. van den Bergh \cite{vdb} (and the authors in \cite{Willwacher2}) that $w_{j}=-(-1)^{j(j-1)/2}\frac{B_{j}}{2j\, j!}$.\footnote{This differs slightly from the formula in \cite{vdb} since we included the symmetry factor $j$ of the wheel in the weight.}
Hence it follows that
\[
[\tilde{\mV}_0(1)] = [e^{-\sum_{j\geq 1} \frac{B_{2j}}{4j\, (2j)!}u^{2j}\tr(R^{2j})}] = 
\sideset{}{^\frac{1}{2}}\det\left(\frac{u R/2}{\sinh(u R/2)} \right)
= \hat{A}_u(M)\, .
\]
\end{proof}

\section{Compatibility of the cyclic Shoikhet--Dolgushev map with cap products}
\label{s-4}
The goal of this section is to define the morphisms $\mV_I^\pi$ and $\mV_H^\pi$ occuring in Proposition \ref{prop2} above and give a proof of that proposition. 
Both morphisms will be defined as a sum of terms corresponding to Shoikhet graphs. In fact the formulas will be identical to those for the twisted version $\mV^\pi$ of B. Shoikhet's morphism $\mV$ \cite{Shoikhet}, except that the weights are computed slightly differently. Recall that the Taylor component $\mV_1^\pi$ of this morphism is defined by the formula
\begin{multline*}
 \mV_1^\pi(\gamma;a_0, \dots, a_n)
 = \\ =
 \sum_{m\geq 0} \frac{\hbar^m}{m!}
\sum_{\Gamma\in G(m+1,n)}
\left(\int_{\pi_1^{-1}(U_{1})\subset C_{m+1,n} }
\omega_\Gamma
\right)
D_\Gamma(\gamma,\pi,\dots, \pi;a_0,\dots, a_n)
\, . 
\end{multline*}
Here $\gamma\in \Tpoly$ and $a_0,\dots, a_n$ are smooth functions. The sum is over all Shoikhet graphs, see Appendix \ref{app:graphs} or \cite{Shoikhet}. The form valued multi-differential operator $D_\Gamma$ is defined in equation \eqref{equ:dgammadef}. The integral is over a certain subspace $\pi_1^{-1}(U_{1})$ of the space $C_{m+1,n}$ (see Definition \ref{def:configspace}), which is essentially B. Shoikhet's space of configurations of $m+n+3$ points in the unit disk. Here ``essentially'' means that our space $C_{m+1,n}$ is bigger by a factor $S^1$, which we interpret as a choice of framing (direction) at the center of the disk. The map $\pi_1: C_{m+1,n}\to C_{1,n}$ is the forgetful map, forgetting the positions of $m$ points $2,\dots,m+1$. The subspace $U_1$ is obtained by fixing the framing to align with the position of point $\bar 0$ on the boundary of the disk. For the notation used, see Appendix \ref{app:graphs}, where also the differential form $\omega_\Gamma$ is defined. 

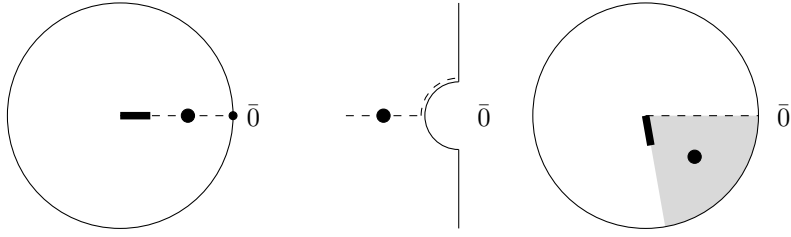
\begin{figure}
\begin{center}
%
\begin{tikzpicture}[scale=.5]
 \draw (0,0) circle (3);
 \draw[dashed] (0,0)--(0:3);
 \draw[line width=1mm] (0,0)--(0:.8);
 \node[circle, draw, fill, minimum size=5, inner sep=0] at (0:1.8) {};
 \node[label=0:{$\bar 0$}, circle, draw, fill, minimum size=3, inner sep=0] at (0:3) {};

\draw(9,3) -- (9,.9);
\draw(9,-.9) -- (9,-3);

 \node[label=0:{$\bar 0$}] at (9,0) {};
\draw[dashed](6,0) -- (8,0);
\draw[dashed](8,0) arc (180:90:1);
\draw(9,-.9) arc (270:90:.9);
 \node[circle, draw, fill, minimum size=5, inner sep=0] at (0:7) {};
\end{tikzpicture}
\hspace{.2cm}
\begin{tikzpicture}[scale=.5, yscale=-1]
\path[fill=gray!30] (0,0) -- ++(0:3) arc (0:80:3) -- cycle;
\draw (0,0) circle (3);
\draw[dashed] (0,0)--(0:3);
 \draw[line width=1mm] (0,0)--(80:.8);
 \node[circle, draw, fill, minimum size=5, inner sep=0] at (40:1.7) {};
 \node[label=0:{$\bar 0$}] at (0:3) {};
\end{tikzpicture}
\end{center}
\caption{\label{fig:subspaces} The subspaces $U_I$ (left) and $U_H$ (right). On the left, the inner and outer framings are aligned and the vertex is restricted to this common line. 
The middle picture shows a magnification of the left in an infinitesimal neighborhood of $\bar 0$.
On the right, the framings are not aligned and the vertex may move in the shaded region between them. Note that on the right the framings are not fixed, but part of the configuration.}
\end{figure}

The formulas for $\mV_I^\pi$ and $\mV_H^\pi$ are given by almost the same expression, except that the subspace $U_{1}$ is replace by different subspaces. Concretely
\begin{multline*}
 \mV_{\mathit{Op}}^\pi(\gamma;a_0, \dots, a_n)
 = \\ =
 \sum_{m\geq 0} \frac{\hbar^m}{m!}
\sum_{\Gamma\in G(m+1,n)}
\left(\int_{\pi_1^{-1}(U_{\mathit{Op}})\subset C_{m+1,n} }
\omega_\Gamma
\right)
D_\Gamma(\gamma,\pi,\dots, \pi;a_0,\dots, a_n)
\end{multline*}
where $\mathit{Op}$ stands for either $I$ or $H$. The subspaces $U_I$ and $U_H$ corresponding to the operation $I$ and $H$ are depicted in Figure \ref{fig:subspaces}. 
For $U_I\subset U_1$ one fixes the framing to the position of point $\bar 0$ on the boundary, and additionally restricts the point $1$ in the interior to that line segment, connecting the center point $0$ with $\bar 0$. For $U_H$, there is no constraint on the inner framing, but the point $1$ is restricted to the disk segment between (in clockwise sense) $\bar{0}$ and the framing at $0$. This is the region colored gray in Figure \ref{fig:subspaces}. The operation $\mV_I$ has already been found by D. Calaque and C. Rossi \cite{CalaqueRossi}, cf. the Remark before section \ref{sec:imof1}.
%
The main result of this section is the following theorem.
\begin{thm}
\label{thm:compat}
The operation $\mV_H^\pi(\gamma)$ satisfies
\begin{multline*}
 \mV_0^\pi H_{\mU_1^\pi(\gamma)} 
- H_\gamma \mV^\pi_0 (-1)^{deg_H}
- \mV^\pi_1(\gamma)
+ (-1)^{deg_f} d \mV_I^\pi(\gamma)
+\mV_I^\pi(\gamma)B (-1)^{deg_H}=
\\=
-(-1)^{|\gamma|} \mV^\pi_H(\gamma) b_\star
- L_\pi \mV^\pi_H(\gamma)
-\mV^\pi_H(\co{\pi}{\gamma})
 \end{multline*}
where $deg_H$ and $deg_f$ are the Hochschild and differential form grading operators.
\end{thm}
\begin{proof}
 This is a standard Kontsevich-Stokes proof. Concretely, one applies Stokes' Theorem to the integrals
\[
 \int_{\pi_1^{-1}(U_H)} d\omega_\Gamma = \int_{\p \pi_1^{-1}(U_H)} \omega_\Gamma\, .
\]
The left hand side vanishes since the weight forms are closed. Computing the right hand side yields some polynomial relations between weights, which are equivalent to the equation in the Theorem. 
In our case, the following codimension 1 boundary strata appear.
\begin{enumerate}
\item $k\geq 2$ interior vertices can approach each other in the interior. By a Kontsevich lemma, the contribution is zero unless $k=2$. If the vertex 1 (corresponding to $\gamma$) is one of the two vertices collapsing, the result is the term $\mV^\pi_H(\co{\pi}{\gamma})$. If not, the result is zero by the Jacobi identity $\co{\pi}{\pi}=0$.
\item Some interior vertices (except 1) approach the center 0. This produces the term $L_\pi \mV^\pi_H(\gamma)$, see \cite{Shoikhet}.
\item \label{item:diffbdry} A cluster of $k$ interior vertices, one of them is 1, approaches the center 0. If $k>1$ this yields zero contribution by (a slight variation of) Kontsevich's Lemma. If $k=1$, i.e., only the vertex 1 approaches 0, we claim that the contribution is $\pm H_\gamma \mV^\pi_0$. One might take this as the definition of $H_\gamma$. But let us nevertheless show that the term produced equals $\pm\frac{1}{2}d L_\gamma$. Let $\alpha$ be the angle between the inner and outer framing. Start by integrating out the position of vertex 1. As in the previous stratum, this produces a term $L_\gamma$, but there remains a factor $\alpha$ in the weight form due to the restricted integration domain. Integrating out the inner framing yields an overall de Rham derivative (as in \cite{Willwacher}), weighted by a factor $\frac{1}{2}=\int_0^1 \alpha d\alpha$.
\item The vertex 1 (corresponding to $\gamma$) can approach the line connecting $0$ and $\bar{0}$. In this case one can easily integrate out the framing at the center vertex (see \cite{Willwacher}) and obtains the term
$d \mV_I^\pi(\gamma)$.
\item The vertex 1 can approach the line corrsponding to the framing at $0$. The result is equal to
$\mV_I^\pi(\gamma)B$.
\item The inner and outer framing align. This leaves no restriction on the position of the vertex 1 and hence produces $\mV^\pi_1(\gamma)$.
\item Vertices on the boundary circle and zero or more interior vertices approach each other. If the vertex 1 (corresponding to $\gamma$) is not in the approaching cluster, this produces the term $\mV^\pi_H(\gamma) b_\star$.
\item If the vertex 1 is in the approaching cluster, the resulting term is equal to $\mV_0^\pi H_{\mU_1^\pi(\gamma)}$.
\end{enumerate}
We apologize for not being very precise about signs.
\end{proof}

To globalize the above result, one has to replace multivector fields, multidifferential operators, differential forms and Hochschild chains by their Fedosov resolutions and apply the above map fiberwise, as usual. See, e.g., \cite{Dolgushev} for a thorough account of the standard globalization procedure.

\subsection{(Remainder of) Proof of part (2) of Proposition \ref{prop2}}
The above result shows that the equality in Proposition \ref{prop2} is satisfied in linear order in $u$. For the sign, note that in our case $\gamma=\dot{\pi}$ is even, as are all forms and Hochschild chains involved. 

The quadratic order in $u$ of the equality in Proposition \ref{prop2} asserts that the homotopy $\mV_H^\pi$ is again compatible with de Rham and Rinehart differential, i.e. that
\[
 d\mV_H^\pi(\gamma) + \mV_H^\pi(\gamma)B=0\, .
\]
Actually, both terms vanish separately.\footnote{This could be derived form the $KS$ relations, together with the facts $B\circ H=H\circ B=0$ and the vanishing of the homotopy corresponding to the composition. However, we do it by hand here.}
We restrict ourselves to the local case: Integrating out the inner framing in the explicit integral formula for $\mV_H^\pi(\gamma)$ shows that $\mV_H^\pi(\gamma)=d(\dots)$ and hence the first term vanishes. On the other hand, consider $\mV_H^\pi(\gamma)B$. It consists of terms having $1$ inserted at the marked point on the boundary of the disk (the outer framing). Hence there is no angle form depending on the outer framing and the integral vanishes. \hfill \qed

\subsection{Globalization and the proof of part (1) of Proposition \ref{prop2}}
Let us next turn to the proof of part (2) of Proposition \ref{prop2}. By the previous arguments, we know that part (2) of that proposition holds. Hence we already know that part (1) holds in the local case, i.e., $M=\R^n$. What we need to check is that the statement survives globalization. 
Let us begin by introducing some notation, which will simplify the following discussion. 

\begin{defi}
 Let $\alg{g}$ be a dg Lie, or $L_\infty$ algebra and let $\alg{m}$ be a dg Lie, or $L_\infty$ module over $\C((u))$. We call $\alg{m}$ a \emph{GM module}\footnote{GM here stands for Gauss-Manin. We made up the name just to convey the correct intuition in the present context, it is not historically motivated.} if it is additionally endowed with an operation 
\begin{align*}
 \hat I \colon \alg{g}\otimes\alg{m} \to \alg{m}
\end{align*}
such that 
\begin{equation}
\label{equ:DhatI}
 \co{D}{\hat I_\gamma} + \hat I_{b \gamma} = uL_\gamma 
\end{equation}
where $D$ is the differential on $\alg{m}$, $L$ denotes the action of $\alg{g}$ on $\alg{m}$ and $b$ is the differential on $\alg{g}$. 
\end{defi}
Let $\alg{m}'$ be another GM module and suppose $\mV:\alg{m}\to\alg{m}'$ is an $L_\infty$ morphism. 
It is then easy to check that for any cocycle $\gamma\in \alg{g}$, the operation
\[
 \Delta_\gamma := \co{\hat I_\gamma}{\mV_0} + u\mV_1(\gamma)
\]
commutes with the differential. Hence it induces a map on homologies $H(\alg{m})\to H(\alg{m'})$. 

\begin{defi}
In the situation above, we say that $\mV$ is \emph{compatible with the Gauss-Manin connection} if the operations $\Delta_\gamma$ vanish on homology, for all cocycles $\gamma\in \alg{g}$.
\end{defi}

The differential forms $\Omega((u))$ and the the cyclic chains $PC$ are both GM modules, with operations $\hat I$ as introduced in the previous sections.
Our goal is to show that the map $\mV^\pi: PC[[\hbar]] \to \Omega((u))[[\hbar]]$ is compatible with the Gauss-Manin connection.
Let us make some elementary observations we will need later.

\begin{lemma}
\label{lem:compositionGM}
Let $\alg{g}$ be a Lie or $L_\infty$ algebra and let $\alg{m}$, $\alg{m}'$, $\alg{m}''$ be GM modules. Let $f:\alg{m}\to\alg{m}'$ and $g:\alg{m}'\to \alg{m}''$ be ($L_\infty$) morphisms of modules, that are compatible with the Gauss-Manin connection. Then $g\circ f: \alg{m}\to\alg{m}''$ is also compatible with the Gauss-Manin connection.
\end{lemma}
\begin{lemma}
\label{lem:pullbackGM}
Let $\alg{g}$ be a Lie or $L_\infty$ algebra, let $\alg{m}$, $\alg{m}'$ be GM modules, and let $f:\alg{m}\to\alg{m}'$ be an $L_\infty$ morphism compatible with the Gauss-Manin connection. Let $\alg{g'}$ be another $L_\infty$ algebra and $F: \alg{g'}\to\alg{g}$ be an $L_\infty$ morphism. Than the pulled back $L_\infty$ modules $F^*\alg{m}$, $F^*\alg{m}'$ over $\alg{g}'$ carry natural structures of GM modules, and the map $F^*f: F^*\alg{m}\to F^*\alg{m}'$ is compatible with the Gauss-Manin connection.
\end{lemma}

\begin{rem}
 The notations ``GM module'' and compatibility with the ``Gauss-Manin connection'' used here should not be taken too serious. 
In all cases we apply the above lemmas, $\hat I$ will be the connection form of the Gauss-Manin connection, but in general, it is actually not justified to use this name.
\end{rem}

\subsection{Reminder on the Dolgushev-Fedosov globalization procedure}
\label{sec:globreminder}
Let us briefly recall V. Dolgushev's globalization procedure. For details, we refer to \cite{Dolgushev}. 
Fix an affine, torsion-free connection on $M$.
Let $\Tpoly^{\formal}$, $\Omega^\formal$, $\Dpoly^\formal$ and $C^\formal$ be V. Dolgushev's resolutions of $\Tpoly$, $\Omega$, $\Dpoly$ and $C=C(A)$.\footnote{In this subsection we will omit the symbol ``${}^\bullet$'' from complexes to avoid cluttering the notation too much. E.g., we write $\Omega$ instead of $\Omega^{-\bullet}$.} The differentials on these spaces depend on the choice of connection. We will denote the part of the differential coming from the connection by $\nabla$, though this is slightly misleading.
One has the following diagram
\[
\begin{CD}
(\Tpoly,0) @>F>> (\Tpoly^\formal, \nabla) @>>> (\Dpoly^\formal, \nabla+b) @<G<< (\Dpoly, b) \\
\Downarrow @. \Downarrow @. \Downarrow @. \Downarrow \\
(\Omega,0) @>f>> (\Omega^\formal, \nabla) @<\mV_f<< (C^\formal, \nabla+b) @<g<< (C,b).
\end{CD}
\]
Here the objects in the upper row are dg Lie algebras, connected by a chain of $L_\infty$ quasi-isomorphisms. The left- and rightmost maps ($F$, and $G$) are the ``flat lifts'', and are ordinary dg Lie algebra morphisms. The middle map is obtained by applying M. Kontsevich's $L_\infty$-morphism fiberwise. The objects in the bottom row are dg Lie algebra modules, the double arrows shall symbolize the action. The arrows in the bottom row are $L_\infty$ morphisms of modules. For example, the left-most arrow in the bottom row denotes an $L_\infty$-morphism of $L_\infty$-modules over $\Tpoly$, where $\Omega^\formal$ is considered $\Tpoly$-module by pulling back the $\Tpoly^\formal$-module structure along the map $\Tpoly\to \Tpoly^\formal$. Again the left- and rightmost maps ($f$ and $g$) in the bottom row are the ``flat lifts'' and contain only a single non-vanishing Taylor component. The middle map is constructed by applying B. Shoikhet's morphism of $L_\infty$-modules fiberwise. From this diagram it follows that there is an an $L_\infty$-morphism $\Tpoly\to \Dpoly$ and an $L_\infty$ morphism of modules $C\to \Omega$. However, to obtain it, one needs to invert the arrows $f$ and $G$, up to homotopy.

Next let us recall the results of \cite{Willwacher}. It is shown there that cyclic differential $B$ and the de Rham differential $d$ are intertwined by Shoikhet's $L_\infty$-morphism. Furthermore, the ``flat lifts'' also commute with these operators. Hence one has the following diagram of $L_\infty$ quasi-isomorphisms
%

\[
 \begin{CD}
 (\Tpoly,0)  @>F>> (\Tpoly^\formal, \nabla) @>>> (\Dpoly^\formal, \nabla+b) @<G<< (\Dpoly, b) \\
\Downarrow @. \Downarrow @. \Downarrow @. \Downarrow \\
(\Omega((u)),u d) @>f>> \substack{ (\Omega^\formal((u)), \\ \nabla+u d_f)} @<<< \substack{ (PC^\formal, \\ \nabla+ b+ uB_f)} @<g<< (PC,b+ u B).
 \end{CD}
\]
Here $PC$ ($PC^\formal$) denotes the (resolution of the) periodic cyclic complex. The subscripts ``f'' in $d_f$ and $B_f$ shall indicate that the operators are applied fiberwise. Next one needs to invert the arrows $f$ and $G$. By abstract reasons, there is an $L_\infty$ quasi-isomorphism $\tilde f: \Omega^\formal((u))\to \Omega((u))$ such that $\tilde f \circ f$ is homotopic to the identity. For us, this means that there is an $L_\infty$ quasi-isomorphism 
\[
 \Phi \colon \Omega((u)) \to \Omega((u))[t,dt]
\]
such that $\Phi(t=0)=\mathit{id}$ and $\Phi(t=1) = \tilde f \circ f$.
Similarly, the morphism $G$ can be inverted up to homotopy to an $L_\infty$ quasi-isomorphism $\tilde G$, such that $G\circ \tilde G$ is homotopic to the identity. However, we need to simultaneously modify the morphism on modules $g$. Note that the naive guess, i.e., the pullback of $g$ by $\tilde G$, write $\tilde G^* g$, will only produce a morphism of modules 
\[
 PC \to (G\circ \tilde G)^* PC^\formal\,.
\]
Here the notation means that the $L_\infty$-module structure on $PC^\formal$ is changed by pull-back along the endomorphism $G\circ \tilde G$ of $\Dpoly^\formal$.
The $L_\infty$ modules $(G\circ \tilde G)^* PC^\formal$ and $PC^\formal$ are $L_\infty$-isomorphic, but one should be careful and choose the isomorphism as follows: The structures of an $L_\infty$ algebra $\alg{g}$ and an $L_\infty$ module $\alg{m}$ can be packaged into a single $L_\infty$ structure on the semi-direct product $\alg{g}\ltimes\alg{m}$. In our case, $\alg{g}=\Dpoly^\formal$ and $\alg{m}=PC^\formal$. The homotopy linking $G\circ \tilde G$ to the identity can be extended to an $L_\infty$ quasi-isomorphism 
\[
 \Psi \colon \alg{g}\ltimes\alg{m} \to \alg{g}\ltimes\alg{m}[t, dt]
\]
such that $\Psi(t=0)$ is the identity and the $\alg{g}$-component of $\Psi(t=1)$ is $G\circ \tilde G$. 
Put an auxiliary grading on $\alg{g}\ltimes\alg{m}$ by declaring $\alg{m}$ to be concentrated in degree $1$ and $\alg{g}$ in degree $0$.
Then $\Psi$ can be chosen to respect this grading. In other words its Taylor components vanish if more than one input is in $\alg{m}$ and have output in $\alg{m}$ iff exactly one input is in $\alg{m}$.

The $\alg{m}$-component of $\Psi(t=1)$ gives an $L_\infty$ isomorphism 
\[
PC^\formal \to (G\circ \tilde G)^* PC^\formal
\]
which can be inverted and composed with $(G\circ \tilde G)^*g$ to yield the desired morphism, say $\tilde g$.
The $L_\infty$ morphism 
\[
 \mV \colon (PC,b+ u B) \to (\Omega((u)),u d)
\]
we use in the present paper is then the composition of the quasi-isomorphisms
\[
 \mV = \tilde f \circ \mV_f \circ \tilde g.
\]
Note that there is no canonical choice for the inverses up to homotopy we had to pick above.
Next we pick a formal Poisson structure $\pi \in \hbar \Tpoly[[\hbar]]$ and twist all occuring $L_\infty$ algebras and modules. This yields the following diagram

\begin{equation}
\label{equ:dolgushevdiagram}
 \begin{CD}
  (\Tpoly[[\hbar]],\co{\pi}{\cdot}) @>F>> \substack{ (\Tpoly^\formal[[\hbar]], \\ \nabla+\co{\pi}{\cdot})} @>\mU^\pi>> \substack{(\Dpoly^\formal[[\hbar]], \\ \nabla+b_\star) }@>\tilde G^\pi>> (\Dpoly[[\hbar]], b_\star) \\
\Downarrow @. \Downarrow @. \Downarrow @. \Downarrow \\
\substack{ (\Omega((u))[[\hbar]], \\ L_\pi + u d) } @<{\tilde f^\pi}<<  \substack{(\Omega^\formal((u))[[\hbar]], \\ \nabla+L_\pi + u d_f)} @<{\mV^\pi_f}<< \substack{(PC^\formal[[\hbar]], \\ \nabla+b_\star+ uB_f)} @<\tilde g^\pi<< \substack{ (PC[[\hbar]], \\ b_\star+ u B)}.
 \end{CD}
\end{equation}

The composition of the lower arrows is $\mV^\pi$. Our goal is to show that this morphism $\mV^\pi$ is compatible with the Gauss-Manin connection. By Lemmas \ref{lem:compositionGM} and \ref{lem:pullbackGM} it is sufficient to check this statement for each of the three squares separately, i.e., show that $\tilde f^\pi$, $\mV^\pi_f$ and $\tilde g^\pi$ are each compatible with the Gauss-Manin connection. This is done as follows:

\begin{enumerate}
 \item Note that $\tilde f^\pi$ is the inverse, up to homotopy, of the twisted version $f^\pi$ of the morphism $f$. (In fact, the morphisms $f$ and $f^\pi$ are identitical since $f$, the flat lift, has only a single non-vanishing Taylor component.) It is easy to see that $f^\pi$ is compatible with the Gauss-Manin connection. Hence it suffices to show that compatibility with the Gauss-Manin connection is preserved under taking inverses up to homtopy. This statement is Corollary \ref{cor:inverseGM} of Section \ref{sec:GMhom}.

\item For the middle square, note that $\mV^\pi_f$ is obtained by applying B. Shoikhet's local $L_\infty$ morphism fiberwise. Part (2) of Proposition \ref{prop2} states that this local morphism is compatible with the Gauss-Manin connection, and even gives an explicit homotopy. Using the fiberwise version of this homotopy, one shows that $\mV^\pi_f$ is indeed compatible with the Gauss-Manin connection.

\item For the rightmost square, note that $\tilde g^\pi$ is the composition
\[
 (G^\pi)^* PC[[\hbar]] \to (G^\pi\circ \tilde G^\pi)^* PC^\formal[[\hbar]] \to PC^\formal[[\hbar]]. 
\]
The left hand morphism is a pullback of the flat lift, which is compatible with the Gauss-Manin connection on the nose. Hence, by Lemma \ref{lem:pullbackGM} it is itself compatible with the Gauss-Manin connection. It will be shown in section \ref{sec:alghomGM} that the right hand morphism as well is compatible with the Gauss-Manin connection.
\end{enumerate}

\subsection{A conceptual remark, and the Kontsevich-Soibelman operad}
\begin{figure}
\begin{center}
\psset{unit=.8cm}
\psset{arrowscale=1}
\begin{pspicture}(-2,-2.25)(2,2.25) 
\psellipse(0,-1.75)(2,.5)
\psellipse(0,1.75)(2,.5)
\psline(-2,1.75)(-2,-1.75)
\psline(2,1.75)(2,-1.75)

\cnode*(0,1.25){3pt}{ain}
\cnode*(0,-2.25){3pt}{aout}
\cnode(1,0){3pt}{vin}
\cnode(0,-0.5){1pt}{mult}

\ncline{ain}{aout}
\ncline{vin}{mult}
\uput{1ex}[ 90](ain){in}
\uput{1ex}[-90](aout){out}
\end{pspicture}
\hspace{.2cm}
\psset{unit=.8cm}
\psset{arrowscale=1}
\begin{pspicture}(-2,-2.25)(2,2.25) 
\psellipse(0,-1.75)(2,.5)
\psellipse(0,1.75)(2,.5)
\psline(-2,1.75)(-2,-1.75)
\psline(2,1.75)(2,-1.75)

\cnode*(0.5,1.266){3pt}{ain}
\cnode*(0,-2.25){3pt}{aout}
\cnode(1,0){3pt}{vin}
\cnode*(0,-1){3pt}{theone}

\cnode(1,-2.183){1pt}{vcon}
\cnode(0.5,-2.234){1pt}{incon}

\ncline{vin}{vcon}
\ncline{ain}{incon}
\ncline{theone}{aout}
\uput{1ex}[ 90](ain){in}
\uput{1ex}[-90](aout){out}
\uput{1ex}[ 90](theone){\bf 1}
\end{pspicture}
\end{center}
\caption{\label{fig:KSops} The operations $I$ (cap product, left) and $H$ (right) in the Kontsevich Soibelman operad.}
\end{figure}
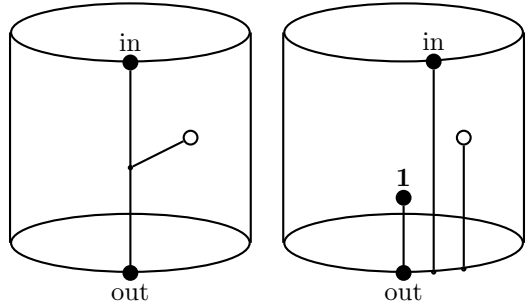

The above proof is conceptually not very satisfying. The clean way to do this proof is to show that (i) M. Kontsevich's morphism and the Shoikhet-Dolgushev morphism can be extended to a morphism of algebras over the Kontsevich-Soibelman (KS) operad \cite{KontsevichS2}, up to homotopy, and (ii) that the $L_\infty$ part of a KS map up to homotopy is compatible with the Gauss-Manin connection. Operations in the KS operad are given by certain graphs. For example, the graphs corresponding to the operations $I$ and $H$ above are depicted in Figure \ref{fig:KSops}. To show the existence of a KS map up to homotopy, one has to construct a tower of homotopies, for the various operations in the KS operad. We construct here only a low order homotopy, corresponding to the operation $H$.
The full construction will be published in a forthcoming paper.

\subsection{Compatibility of the Gauss-Manin connection with $L_\infty$-homotopies}
\label{sec:GMhom}

\begin{lemma}
 Suppose $\alg{g}$ is an $L_\infty$ algebra and $\alg{m}$, $\alg{m'}$ are GM modules. Suppose $\mV,\mW : \alg{m}\to\alg{m}'$ are two homotopic $L_\infty$ morphisms of modules. Then $\mV$ is compatible with the Gauss-Manin connection iff $\mW$ is.
\end{lemma}
\begin{rem}
 For us, $\mV,\mW$ are homotopic if there is an $L_\infty$ morphism $f: \alg{m}\to\alg{m}'[t, dt]$ such that $f(t=0)=\mV$, $f(t=1)=\mW$.
\end{rem}
\begin{proof}
 Let $f$ be as in the remark. It fits into a diagram, commutative in homology
\[
 \alg{m} \stackrel{f}{\to} \alg{m}'[t, dt] \rightrightarrows \alg{m}'
\]
where the two right hand arrows are the evaluations at $t=0$ and $t=1$ respectively. $\alg{m}'[t, dt]$ is a GM module by $t$ linear extension of the GM structure on $\alg{m'}$. Let us denote the operator giving the GM structure also by $\hat I$. We know the for a closed $\gamma\in \alg{g}$ the morphism 
\[
 \Delta_\gamma := \co{\hat I_\gamma}{f_0} + uf_1(\gamma)
\]
commutes with the differential and hence defines an operation $[\Delta_\gamma]$ on homology. It fits into a commutative diagram 
\[
 H(\alg{m}) \stackrel{[\Delta_\gamma]}{\to} H(\alg{m}'[t, dt]) \rightrightarrows H(\alg{m}').
\]
Hence the lower composition is zero iff the upper composition is. This proves the Lemma.

\end{proof}

\begin{cor}
\label{cor:inverseGM}
 Suppose $\alg{g}$ is an $L_\infty$ algebra, $\alg{m}$, $\alg{m'}$ are GM modules and $\mV : \alg{m}\to\alg{m}'$ is an $L_\infty$ quasi-isomorphism, compatible wih the Gauss-Manin connection. Then any inverse up to homotopy 
$\tilde \mV : \alg{m}'\to\alg{m}$
is also compatible with the Gauss-Manin connection.
\end{cor}
\begin{proof}
 We have to show that 
\[
 \Delta_\gamma := \co{\hat I_\gamma}{\tilde \mV_0} + u\tilde \mV_1(\gamma)
\]
vanishes on homology. Since $\mV$ is a quasi-isomorphism, this is equivalent to saying that 
\[
 \left(\co{\hat I_\gamma}{\tilde \mV_0} + u\tilde \mV_1(\gamma)\right) \circ \mV_0
=
\co{\hat I_\gamma}{\tilde \mV_0\circ \mV_0} + u\tilde \mV_1(\gamma)\circ \mV_0 + u\tilde \mV_0\circ\mV_1(\gamma)
+\co{D}{\cdots}
\]
vanishes on homology. For the equality we used that $\mV$ is compatible with the Gauss-Manin connection. But vanishing of the right hand expression on homology just means that the composition $\tilde \mV\circ \mV$ is compatible with the Gauss-Manin connection. But this follows from the Lemma since $\tilde \mV\circ \mV$ is homotopic to the identity by definition of ``inverse up to homotopy''.
\end{proof}

This Corollary in particular asserts that the map $\tilde F^\pi$ in the left hand square of the diagram \eqref{equ:dolgushevdiagram} is compatible with the Gauss-Manin connection.

\subsection{Compatibility of the Gauss-Manin connection with \emph{algebra} homotopies}
\label{sec:alghomGM}
Let again $\alg{g}$ be an $L_\infty$ algebra, and $\alg{m}$ a GM module. Suppose we have an $L_\infty$ endomorphism $\Psi: \alg{g}\ltimes \alg{m} \to \alg{g}\ltimes \alg{m}$ homotopic to the identity. We also assume, as usual, that $\Psi$ and the homotopy respect the $\alg{m}$ degree, i.e., all Taylor components vanish when at least two of its arguments are in $\alg{m}$, and the output is $\alg{m}$ if exactly one input is.
In this case $\Psi$ encodes an $L_\infty$ morphism (of $L_\infty$ algebras) $F^1: \alg{g}\to \alg{g}$, together with a morphism of $L_\infty$ modules $f^1: \alg{m}\to (F^1)^*\alg{m}$.

\begin{lemma}
  In the situation above, $f^1$ is compatible with the Gauss-Manin connection.
\end{lemma}
\begin{proof}
 Since $\Psi$ is homotopic to the identity, we can extend $\Psi$ to a morphism $\tilde \Psi\colon \alg{g}\ltimes \alg{m} \to \alg{g}\ltimes \alg{m}[t,dt]$, such that the compositions with the evaluation at $t=0$ and $t=1$ are the identity and $\Psi$ respectively. In particular there is an $L_\infty$ morphism $F: \alg{g}\to \alg{g}[t,dt]$, and an $L_\infty$ morphism of $L_\infty$ $\alg{g}$-modules $f: \alg{m}\to F^*\alg{m}[t,dt]$, such that the evaluations at $t=0$ and $t=1$ are the identities and $F^1$, $f^1$ respectively. Here, as before, $F^*\alg{m}[t,dt]$ is $\alg{m}[t,dt]$, endowed with a $\alg{g}$-module structure by pulling back the natural $\alg{g}[t,dt]$ module structure via $F$. Note that $F^*\alg{g}[t,dt]$ is also a GM module, with the pulled back GM structure, and that the evaluation maps are compatible with the Gauss-Manin connection. Denote the GM structure on $\alg{m}[t,dt]$ again by $\hat I$, and for a closed $\gamma \in \alg{g}$ define as before 
\[
 \Delta_\gamma := \co{\hat I_\gamma}{f_0} + uf_1(\gamma).
\]
It induces a map $[\Delta_\gamma]$ on homology, whose composition with the evaluation at $t=0$, i.e., 
\[
  H(\alg{m}) \stackrel{[ \Delta_\gamma ]}{\longrightarrow} H(\alg{m}[t, dt]) \stackrel{\mathit{ev}_{t=0}}{\longrightarrow} H(\alg{m}').
\]
is zero. Since the evaluations are quasi-isomorphisms it follows that $[\Delta]=0$, i.e., $f$ is compatible with the Gauss-Manin connection. But then also $f^1 = \mathit{ev}_{t=1} \circ f$ is compatible with the Gauss-Manin connection.
\end{proof}

A twisted version of this Lemma produces the main result of this section, finishing the proof of Proposition \ref{prop2}.

\begin{lemma}
 The map $\tilde g^\pi$ of section \ref{sec:globreminder} is compatible with the Gauss-Manin connection.
\end{lemma}
\begin{proof}
 Using the notation of Section \ref{sec:globreminder}, the map $\tilde g^\pi$ is a composition
\[
 (G^\pi)^* PC[[\hbar]] \to (G^\pi\circ \tilde G^\pi)^* PC^\formal[[\hbar]] \to PC^\formal[[\hbar]]. 
\]
The left hand morphism is a pullback of the flat lift $g$, which is compatible with the Gauss-Manin connection on the nose. Hence, by Lemma \ref{lem:pullbackGM} it is itself compatible with the Gauss-Manin connection. To show that the right hand morphism is compatible, it is sufficient to show that its homotopy inverse 
\[
 PC^\formal[[\hbar]] \to (F^{1,\pi})^* PC^\formal[[\hbar]]
\]
is. Here $F^{1,\pi}$ is the twisted version of $F^1 = G\circ \tilde G$. By construction of this map in Section \ref{sec:globreminder}, there are quasi-isomorphisms 
\[
PC^\formal[[\hbar]] \to (F^\pi)^* PC^\formal[t,dt][[\hbar]]\rightrightarrows PC^\formal[[\hbar]]
\]
such that the upper composition is the identity, and the lower is $\tilde G^\pi$. Here the two morphisms on the right are the evaluations at 0 and 1. Note that $PC^\formal[[\hbar]]$ is endowed with a GM structure. Hence also $(F^m)^*PC^\formal[t,dt][[\hbar]]$ is a GM algebra by $t$ linear extension and pullback. Since all maps are quasi-isomorphisms and the diagram is commutative in homology, the statement of the Lemma then follows in exactly the same way as in the previous proof.
\end{proof}

\section{Relation to the Tamarkin
-
Tsygan index Theorem}
As observed by D. Tamarkin and B. Tsygan \cite{TamarkinT} the results on cyclic formality may be used to prove an algebraic index theorem as follows.
Let $\pi$ be a Poisson structure and $\pi_t = t\pi$ for $t\in \R$. Let $A_t$ be the space $A=C^\infty(M)[[\hbar]]$ with the star product corresponding to $\pi$.
Let $c\in PH_\bullet(A_1)$ be a cyclic homology class. By parallel transport wrt. the Gauss\ndash Manin connection (see Theorem \ref{t-1}) one obtains a family of cyclic homology classes $c(t)\in PH_\bullet(A_t)$ such that $c(1)=c$. In particular $c(0)\in PH_\bullet(A)$. By Connes isomorphism (see section \ref{s-3}) $PH_\bullet(A)\cong H^\bullet(M)((u))$. We denote the image of $c(0)$ by $ch(c)\in H^\bullet(M)((u))$. 

On the other hand one has a map from periodic cyclic to cyclic homology $PH_\bullet(A_1)\to CH_\bullet(A_1)$. On chains, this map is given by sending positive powers of $u$ to zero. Furthermore, the zeroth cyclic homology is isomorphic to the zeroth Hochschild homology $CH_0(A_1)\cong HH_0(A_1)$. By the Shoikhet\ndash Dolgushev morphism $\mV^\pi$ the latter is, in turn, isomorphic to zeroth Poisson homology $HP_0(M)$. Assume now that there is a volume form $\Omega$ on $M$ and $\pi$ is a unimodular with $\dv_\Omega\pi=0$.\footnote{A Poisson structure is called unimodular if its divergence with respect to some (and hence every) volume form is a Hamiltonian vector field. If this is the case, rescaling the volume form by the exponential of the Hamiltonian function yields a volume form with respect to which the Poisson structure is divergence-free. If $\{\ ,\ \}$ denotes the Poisson bracket, then the last condition is equivalent to $\int_M\{f,g\}\,\Omega=0$, $\forall f,g$.}
Than there is a natural map $HP_0(M)\to \C$ by $f\mapsto \int_M f\Omega$.

Following Tamarkin and Tsygan, let us denote the composition by\footnote{Note that by form degree reasons, terms with positive powers of $u$, though present in the integrand, do not contribute to $I$.} 
\begin{gather*}
I\colon PH_0(A_1) \to CH_0(A_1) \stackrel{\mV^\pi_0}{\to} HP_0(M) \to \C \\
I(c) = \int_M \mV^\pi_0(c) \Omega
\, .
\end{gather*}

By Theorem \ref{t-1} we know that 
\[
[\mV^\pi_0(c)] = [e^{-\iota_\pi /u}\mV^{\pi=0}_0(c(0))] \, .
\]
Furthermore, by Proposition \ref{prop:SDAhat}
\[
\mV^{\pi=0}_0(c(0)) = \hat{A}_u(M) ch(c)\,.
\]
Hence we arrive at the following Theorem, known as the Tamarkin\ndash Tsygan index Theorem. The statement can be found in \cite{TamarkinT}.
\begin{thm}
Let $M$ be a compact manifold, $\pi\in \hbar\Gamma(M,\wedge^2 TM)[[\hbar]]$ a formal Poisson structure, $\Omega$ a volume form on $M$ with $\dv_\Omega\pi=0$ and $c\in PH_0(A_\hbar)$. Then
\[
I(c) = \int_M  \hat{A}_u(M) ch(c) e^{\iota_\pi /u}\Omega\, .
\]
\end{thm}

\appendix 

\section{Graphs and weights}
\label{app:graphs}
We briefly recall here some standard facts about Kontsevich graphs and their weights.
For more details, see \cite{Kontsevich}.

\begin{defi}
\label{def:kontsgraph}
The set $G(m,n)$, $m,n\in \mathbb{N}_0$ of \emph{Shoikhet graphs} consists of directed graphs $\Gamma$ such that
\begin{enumerate}
\item The vertex set of $\Gamma$ is 
\[
V(\Gamma) = \{0,..,m\}\cup \{\bar{0},..,\bar{n}\}
\]
where the vertices $\{1,..,m\}$ will be called the \emph{type I} vertices and the vertices $\{\bar{0},..,\bar{n}\}$ the \emph{type II} vertices.
\item There are no double edges, i.e., edges $(j,w)$ occuring twice in $E(\Gamma)$.
\item Every edge $e=(v,w)\in E(\Gamma)$ starts at a type I vertex, i.e., $v\in\{1,\dots,m\}$.
\item No edge ends at vertex $0$.
\item For each type I vertex $j$, there is an ordering given on
\[
\mathit{Sta}r(j) = \{(j, w) \mid (j, w)\in E(\Gamma),\; w\in V(\Gamma) \}.
\] 
\item There are no tadpoles, i.e., edges of type $(j,j)$.
\end{enumerate}
\end{defi}

Let us next define the (Shoikhet) weight $w_\Gamma$ of $\Gamma\in G(m,n)$. It is an integral of a certain differential form over a compact manifold with corners, the \emph{configuration space} $C_{m,n}$.
\begin{equation}
\label{equ:wgammadef}
w_\Gamma =  \int_{C_\Gamma} \omega_\Gamma
\end{equation}

\begin{defi}
\label{def:configspace}
The {enlarged Shoikhet configuration space} $\tilde{C}_{m,n}$ is the space of embeddings 
\[
(z_0,\dots,z_m,z_{\bar{0}},\dots z_{\bar{n}}): V(\Gamma) \rightarrow D
\]
of the set $\{0,..,m\}\cup \{\bar{0},..,\bar{n}\}$ into the closed unit disk $D=\{z \in \cn{}; |z|\leq 1\}$, together with a distinguished direction at vertex 0, such that
\begin{enumerate}
\item All type I vertices are mapped to the interior of $D$, i.e. $z_j\in D^\circ$ for $j=0,\dots,m$.
\item All type II vertices are mapped to the boundary of $D$, i.e. $z_{\bar{j}}\in \partial D$ for $j=0,\dots,n$.
\item The type II vertices occur in counterclockwise increasing order on the circle, i.e., $0<\arg \frac{z_{\bar{1}}}{z_{\bar{0}}}< \dots < \arg \frac{z_{\bar{n}}}{z_{\bar{0}}}<2\pi$.
\end{enumerate}
The \emph{Shoikhet configuration space} $C_{m,n}$ is the Fulton-MacPherson-Axelrod-Singer compactification of the quotient of $\tilde{C}_{m,n}$ under the action of the automorphism group of the unit disk $\mathit{PSU}(1,1)$. 
\end{defi}
\begin{rem}
 Note that our configuration space is actually 
\[
\{\text{B. Shoikhet's original configuration space}\}\times S^1,
\]
due to the additional data of a framing at $0$. Hence the name is slightly inaccurate.
\end{rem}
\begin{rem}
 Instead of working with equivalence classes of configurations under the $\mathit{PSU}(1,1)$ action, we will generally work with a canonical representative, obtained by fixing the point $z_0$ at the center of the disk, i.e., $z_0=0$, and fixing the point $z_{\bar{0}}$ to be $z_{\bar{0}}=1$.
\end{rem}

The differential form $\omega_\Gamma$ that is integrated over configuration space can be expressed as a product of 
one-forms, one for each edge in $\Gamma$.
\begin{align*}
\omega_\Gamma 
&= \bigwedge_{(0, v)\in E(\Gamma)} \alpha_1(v) \bigwedge_{j=1}^n \bigwedge_{(j, v)\in E(\Gamma)} \alpha(j,v)
\end{align*}

Here the one-form $\alpha_0(v)$ is the differential of the (hyperbolic) angle between the framing at the vertex $z_0$ and the hyperbolic geodesic from $z_0$ to $v$. Similarly, the angle $\alpha(j,v)$ is the differential of the hyperbolic angle between the hyperbolic straight lines $(z_j,v)$ and $(z_j,z_0)$. 

To any graph $\Gamma\in G(m+1,n)$ as above one can associate a function $D_\Gamma$ taking $m$ vector fields $\gamma_1,\dots,\gamma_m$ and $n+1$ functions and returning a $d$-differential form where $d$ is the degree of vertex 0 in $\Gamma$. It is defined such that for a constant $d$-vector field $\gamma_0$
\begin{multline}
\label{equ:dgammadef}
(-1)^d \iota_{\gamma_0} D_\Gamma(\gamma_1\otimes\dots \otimes\gamma_m ; a_1,\dots,a_n)=\\
 = \sum_{\varphi:E(\Gamma)\to [d]} \prod_{j=0}^m (\p_{\varphi(f_1^j)} \p_{\varphi(f_2^j)}\cdots \gamma_j^{\varphi(e_1^j)\varphi(e_2^j)\cdots}) \prod_{k=1}^n (\p_{\varphi(f_1^{\bar{k}})} \p_{\varphi(f_2^{\bar{k}})}\cdots a_k)\, .
\end{multline}
Here the sum runs over all maps $\varphi$ from the edge set of $\Gamma$ to the set $\{1,..,d\}$. The edges incoming to a vertex $v$ are denoted by $f_1^{v},f_2^{v}, \dots$ and the edges outgoing by $e_1^{v},e_2^{v},\dots$

%
\end{document}